%% file: main.tex
\documentclass[hidelinks,12pt]{article}
\usepackage[a4paper, margin=2.5cm]{geometry} 

\usepackage{latexsym,amssymb,graphicx,amsthm,amsmath,mathtools,caption,calc,transparent,float,dsfont,textcomp,comment,bm}
\usepackage[backend=biber,style=alphabetic]{biblatex}
\usepackage[justification= centering]{caption}
\usepackage{hyperref}

\numberwithin{figure}{section}
\numberwithin{equation}{section}

\AtEveryBibitem{
 \clearlist{address}
 \clearfield{date}
 \clearfield{eprint}
 \clearfield{isbn}
 \clearfield{issn}
 \clearlist{location}
 \clearfield{month}
 \clearfield{series}
 \clearfield{doi}
 \clearfield{url}
 
 \ifentrytype{book}{}{
  \clearlist{publisher}
  \clearname{editor}
 }
}

\newtheorem{theorem}{Theorem}[section]

\newtheorem{lemma}[theorem]{Lemma}

\newtheorem{claim}{Claim}[theorem]
\theoremstyle{definition}
\newtheorem{definition}[theorem]{Definition}

\theoremstyle{remark}
\newtheorem{remark}[theorem]{Remark}

\newcommand{\half}{\frac{1}{2}}

\newcommand{\reals}{\mathbb R}

\newcommand{\eps}{\varepsilon}
\newcommand{\ph}{\varphi}
\newcommand{\setin}{\subseteq}
\newcommand{\del}{\partial}

\makeatletter
\newcommand*\twiddle[1]{\mathpalette\wthelper{#1}}
\newcommand*\wthelper[2]{%
        \hbox{\dimen@\accentfontxheight#1%
                \accentfontxheight#11.3\dimen@
                $\m@th#1\widetilde{#2}$%
                \accentfontxheight#1\dimen@
        }%
}

\newcommand*\accentfontxheight[1]{%
        \fontdimen5\ifx#1\displaystyle
                \textfont
        \else\ifx#1\textstyle
                \textfont
        \else\ifx#1\scriptstyle
                \scriptfont
        \else
                \scriptscriptfont
        \fi\fi\fi3
}
\makeatother

\DeclareMathOperator{\lap}{\Delta}

\newcommand{\mfd}{\mathcal{M}}
\newcommand{\heatop}{\left(\del_t-\lap_{g(t)}\right)}

\DeclareMathOperator{\ric}{Ric}

\DeclareMathOperator{\riem}{Rm}

\title{A sharp local maximum principle for Ricci flow}
\author{Laura Flower}
\date{27 November 2026}

\begin{document}

\maketitle

\begin{abstract}
    Localised maximum principles are essential to the study of Ricci flow on noncompact manifolds, with results and applications going back to the work of Perelman. In this work, we obtain a sharp localised maximum principle, improving on a result of Lee--Tam.
\end{abstract}

\section{Introduction}

We say that a smooth family of Riemannian metrics $(g(t))_{t\in[0,T)}$ on a manifold $\mfd$ is a \emph{Ricci flow} if
\[
\del_tg(t)=-2\ric(g(t))
\]
for $t\in[0,T)$. The Ricci flow equation is a nonlinear heat equation for the Riemannian metric, and has seen many successful applications to differential geometry and topology.

It is of interest to study the Ricci flow on noncompact manifolds, where localised maximum principles are essential, with results and applications going back to the work of Perelman \cite{Perelman}. In particular, continuing in the spirit of \cite{BamlerCabezasRivasWilking}, the work of Lee and Tam \cite{LeeTam22} establishes the following localised maximum principle, which we state here in full for the reader's convenience.
\begin{theorem}[{\cite[Theorem 1.1]{LeeTam22}}]\label{PROP Lee Tam thm}
    Let $(\mfd^n,g(t))_{t\in[0,T)}$ be a (possibly incomplete) Ricci flow, with $B_{g(0)}(p,2)\Subset\mfd$ for some $p\in\mfd$. Suppose that on $\mfd\times[0,T)$, we have
    \[
    \ric{g(t)} \leq \alpha t^{-1}
    \]
    for some $\alpha>0$.
    
    Suppose we also have continuous functions $\ph,L:\mfd\times[0,T)\to\reals$,
    with $\ph,L\leq\alpha t^{-1}$,
    such that
    \[
    \heatop\ph\leq L\ph \ \text{ wherever }\ph>0, \text{ in the barrier sense}
    \]
    on $\mfd\times[0,T)$. Finally, suppose that $\ph\leq0$ on $B_{g(0)}(p,2)$ when $t=0$.

    Then for any $\ell>\alpha+1$, there is some time $\widehat{T}=\widehat{T}(n,\alpha,\ell)>0$ such that
    \[
    \ph \leq t^\ell
    \]
    for $t\in[0,T)\cap[0,\widehat{T}]$.
\end{theorem}

Here, and throughout this paper, we use the following definition of barrier sense for differential inequalities.\footnote{This is the same definition as \cite[Definition 18.5]{RFT&A}.}
\begin{definition}\label{DEF barrier sense}
    Given a (possibly nonsmooth) function $f$ and $C\in\reals$, we say that $\heatop f\leq C$ at $(x_0,t_0)$ \emph{in the barrier sense} if the following is true. For any $\eps>0$ there is a $C^2$ function $F$ with $F\leq f$ defined in a neighbourhood of $(x_0,t_0)$, and such that $F(x_0,t_0)=f(x_0,t_0)$, and $\heatop F\leq C+\eps$ at $(x_0,t_0)$.
\end{definition}

Theorem \ref{PROP Lee Tam thm} can in particular be applied to functions related to curvature quantities, and gives the preservation of various notions of non-negative curvature on complete manifolds with scaling-invariant curvature bounds (\cite{LeeTam22}, see also \cites{BamlerCabezasRivasWilking, LiNi20, Wilking2013}). Theorem \ref{PROP Lee Tam thm} has seen numerous further applications to Ricci flow, including obtaining existence results \cites{LeeTopping25Three,LeeTopping25n,Hsiao26,ChauMartens26,Lee23}, uniqueness results and similar \cites{ChanLeePeachey24,Lee26Uniqueness,Flower26,FlowerLeeSchulze26}, and constructing new examples of solitons \cites{ChanLee25Solitons,ChanLeeLai26}. Several of the aforementioned existence results then allow the use of the smoothing properties of Ricci flow to establish geometric results \cites{LeeTopping25Three,LeeTopping25n,Lee23}, and Theorem \ref{PROP Lee Tam thm} has also been used in other cases where geometric results are established using Ricci flow smoothing \cites{ChanLee24,LeeTam25,HuangLiu25,ChanLee25Gap,Lee26Quantification,ChuLeeWan26}.

It is natural to wonder whether a sharp version of Theorem \ref{PROP Lee Tam thm} could be obtained, perhaps by optimising over $\ell$. Simply optimising over $\ell$ in the proof of Lee--Tam does not give the sharp case. However, with more care in the choice of the cutoff functions used and by optimising the bound at various points, in this work we are able to obtain the following.
\begin{theorem}\label{PROP main thm for powers of t}
    There is a universal constant $\widehat{\beta}$ such that the following holds.

    Let $(\mfd^n,g(t))_{t\in[0,T)}$ be a (possibly incomplete) Ricci flow, with $B_{g(0)}(p,r)\Subset\mfd$ for some $p\in\mfd$, $r>0$. Assume that on $\mfd\times[0,T)$, we have
    \[
    \ric{g(t)} \leq \alpha t^{-\gamma}
    \]
    for some $\alpha>0$, $\gamma<2$.

    Suppose we also have continuous functions $\ph,L:\mfd\times[0,T)\to\reals$, and a constant $\beta\leq\widehat{\beta}$, such that
    \[
    L\leq\big(\widehat{\beta}-\beta\big) r^2 t^{-2},
    \]
    \[
    \ph\leq\exp\left(\beta r^2 t^{-1}\right),
    \]
    \[
    \heatop\ph\leq L\ph \ \text{ wherever }\ph>0, \text{ in the barrier sense}
    \]
    on $\mfd\times[0,T)$. Finally, assume $\ph\leq0$ on $B_{g(0)}(p,r)$ when $t=0$.

    Then there exists $\widehat{T}=\widehat{T}(n,\alpha,\gamma,r)$ such that
    \[
    \ph(p) \leq \exp\big(-\big(\widehat{\beta}-\beta\big)r^2t^{-1}\big)
    \]
    for $t\in[0,T)\cap[0,\widehat{T}]$.
\end{theorem}
\begin{remark}
    In particular, $\widehat{T}$ is such that $r=C(\widehat{T}^{1-\gamma/2}+\widehat{T}^{\gamma/2}$) for some $C=C(n,\alpha,\gamma)$. 
\end{remark}
\begin{remark}
    Intuitively, Theorem \ref{PROP main thm for powers of t} says that if $\ph\leq h(t)$ on $B_{g(0)}(p,r)$, and $h(t)$ doesn't grow too quickly as $t\searrow0$, then $\ph(p)\leq h(t)\exp(-Cr^2t^{-1})$. Further, the slower $h(t)$ grows as $t\searrow0$, the quicker the bound on $L$ is allowed to grow as $t\searrow0$.
\end{remark}
\begin{remark}
    Theorem \ref{PROP main thm for powers of t} will be obtained as a corollary of Theorem \ref{PROP main theorem 1/t^2 and exp -1/t}, which allows for a Ricci upper bound of the square of an integrable function of $t$.
\end{remark}

We have therefore replaced the polynomial bound of Theorem \ref{PROP Lee Tam thm}, which gives $\ph\leq t^\ell$ for some $\ell>0$, with the bound $\ph\leq\exp(-Cr^2t^{-1})$ for some $C>0$. This is sharp, as we see from the heat kernel $G:\reals^n\times(0,\infty)\to\reals$ on $\reals^n$, which is given by
\[
G(x,t)=\frac{1}{(4\pi t)^{n/2}}\exp\left(-\frac{|x|^2}{4t}\right),
\]
so
\[
(\del_t-\lap)G=0, \hspace{8ex} \lim_{t\searrow 0}G(\cdot,t)=\delta_0.
\]
If $y\in\reals^n$ satisfies $|y|=r$, then $(\del_t-\lap)G\leq0$ and $G\leq(4\pi t)^{-n/2}\leq\exp(\widehat{\beta}r^2t^{-1})$ in $B(y,r)$ for $t$ sufficiently small, and $G(\cdot,0)=0$ in $B(y,r)$. Then
\[
G(y,t)\geq\exp\left(-\frac{r^2}{4t}\right),
\]
so we cannot do better than a bound which looks like $\exp(-Cr^2t^{-1})$ for some $C>0$.

More generally, the estimate given by Theorem \ref{PROP main thm for powers of t} matches an estimate of Bamler--Cabezas-Rivas--Wilking \cite{BamlerCabezasRivasWilking} for the backwards heat kernel for a Ricci flow with a scaling-invariant curvature bound and scaling-invariant volume lower bound. In particular, let $g(t)_{t\in[0,T)}$ be a complete Ricci flow on a manifold $\mfd$, and let $G(x,t;y,s)$, for $x,y\in\mfd$, $0\leq s<t<T$, be the heat kernel corresponding to the backwards heat equation coupled with the Ricci flow. That is, for fixed $(x,t)\in\mfd\times[0,T)$ we require
\[
(\del_s+\lap_{y,s})G(x,t;\cdot,\cdot)=0, \hspace{8ex} \lim_{s\nearrow t}G(x,t;\cdot,s)=\delta_x.
\]
In particular, for fixed $(y,s)\in\mfd\times[0,T)$, we have
\[
(\del_t-\lap_{x,t}-\operatorname{Scal}(g(t)))G(\cdot,\cdot;y,s)=0, \hspace{8ex} \lim_{t\searrow s}G(\cdot,t;y,s)=\delta_y.
\]
Suppose we have $A,v_0>0$ such that $|\riem(g(t))|\leq At^{-1}$ and $\operatorname{Vol}_{g(t)}(B_{g(t)}(x,\sqrt{t}))\geq v_0t^{n/2}$ for all $(x,t)\in\mfd\times[0,T)$. Then by \cite[Proposition 3.1]{BamlerCabezasRivasWilking}, there is $C(n,A,v_0)>0$ such that for each fixed $y\in\mfd$,
\begin{equation}\label{EQN general heat kernel bound}
    G(x,t;y,0)\leq\frac{C}{t^{n/2}}\exp\left(-\frac{\operatorname{dist}_{g(0)}^2(x,y)}{Ct}\right)
\end{equation}
for every $(x,t)\in\mfd\times[0,T)$. On the other hand, we have
\[
(\del_t-\lap_{x,t})G(\cdot,\cdot;y,0)=\operatorname{Scal}(g(t))G(\cdot,\cdot;y,0)\leq A(n-1)t^{-1}G(\cdot,\cdot;y,0),
\]
and \eqref{EQN general heat kernel bound} gives $G(\cdot,\cdot;y,0)\leq Ct^{-n/2}$ on $\mfd\times[0,T)$. Then the result of Theorem \ref{PROP main thm for powers of t} matches \eqref{EQN general heat kernel bound}.

Finally, it is natural to ask whether the required curvature bound of $\alpha t^{-\gamma}$, for any $\gamma<2$, could be improved to a bound of $\alpha t^{-2}$. However, a curvature bound of $\alpha t^{-2}$ represents a genuine threshold in behaviour, as demonstrated by the Ricci flow coming out of a hyperbolic punctured torus given by Topping \cite{Topping12}. In particular, given a 2-torus with a single puncture and equipped with the hyperbolic metric, as well as the complete self-similarly shrinking solution that one would expect, there is also a Ricci flow $g(t)$ coming out of this space which is compact for all $t>0$. The metric $g(t)$ converges smoothly locally to the metric on the initial space. Further, the flow has curvature decay $\max K \sim t^{-2}$ as $t\searrow0$ \cite{ToppingYin20}. The scalar curvature has evolution equation
\[
\heatop R = 2|\ric|^2 = R^2
\]
since we are in dimension two, and we see that $R$ satisfies a $t^{-2}$ curvature bound and is continuous back to $t=0$ on any compact subset of the initial space, with $R\equiv-1\leq0$ when $t=0$. Therefore, if Theorem \ref{PROP main thm for powers of t} held when we merely had $\ric g(t)\leq \alpha t^{-2}$, then we could apply it in this case on balls with arbitrarily large radius with $\ph=L=R$. Sending the radius of the balls to infinity, we obtain
\[
R\leq0 \hspace{8ex}\text{for all }t\in(0,T).
\]
By Gauss-Bonnet, we therefore have
\[
R=0 \hspace{8ex}\text{for all }t\in(0,T),
\]
contradicting the local smooth convergence of $g(t)$ to the initial metric, since $R\equiv-1$ initially.

The rest of this paper is structured as follows. In \S\ref{SEC cutoff functions}, we construct the cutoff functions which will be needed for the main argument. In \S\ref{SEC existing results}, we briefly recall the method used to prove Theorem \ref{PROP Lee Tam thm} in \cite{LeeTam22}, and optimise over $\ell$ to get a slight improvement. Finally, in \S\ref{SEC main result}, we prove a generalisation of Theorem \ref{PROP main thm for powers of t}. This is achieved for slightly worse bounds on $\ph$ and $L$ by repeating the method in \cite{LeeTam22} using a carefully chosen cutoff function and feeding in the result from \S\ref{SEC existing results}, then carrying out another optimisation step (see Theorem \ref{PROP main theorem 1/t^1+gamma and 1}). The final steps to improve the required bounds on $\ph$ and $L$ makes further use of the continuity of both functions (see Theorem \ref{PROP main theorem 1/t^2 and exp -1/t}). 

\subsection{Acknowledgements}
The author wishes to thank her supervisors, Felix Schulze and Peter Topping,
and Man Chun Lee,
for many helpful and insightful discussions of this topic. The author is supported by the Warwick Mathematics Institute Centre for Doctoral Training, and gratefully acknowledges funding from the University of Warwick and the UK Engineering and Physical Sciences Research Council (Grant number: EP/W524645/1).

\subsection{AI Usage}
No AI was used at any point in the preparation of this paper.

\section{Perelman cutoff functions}\label{SEC cutoff functions}

We begin by discussing the construction of suitable cutoff functions. Our construction is similar to those used in \cite{Perelman}, \cite{Simon13}, \cite{SimonTopping22}, \cite{LeeTam22}, and \cite{LeeTam19}, but we need a refinement of the approach.

The key to the construction is the following lemma of Perelman \cite[Lemma 8.3(a)]{Perelman}.

\begin{lemma}\label{PROP Perelman's lemma}
    Let $(\mfd^n,g(t))_{t\in[0,T)}$ be a (possibly incomplete) Ricci flow. Suppose that at some time $t_0$, we have
    \[
    \ric{g(t_0)} \leq \frac{3}{2}(n-1)K^2 \hspace{3ex}\text{in}\hspace{3ex} B_{g(t_0)}(p,r_0),
    \]
    for some point $p\in\mfd$, some $r_0>0$, and some number $K>0$. Then $d(x,t):=\operatorname{dist}_{g(t)}(p,x)$ satisfies
    \[
    \heatop d \geq -(n-1)\left(K^2r_0+\frac{1}{r_0}\right) \hspace{3ex}\text{outside}\hspace{3ex} B_{g(t_0)}(p,r_0)
    \]
    in the barrier sense.
\end{lemma}

\begin{remark}
    When $r_0=K^{-1}$, we obtain the best possible bound, namely $\heatop d \geq -2(n-1)K$, so for simplicity we always assume this to be the case below. However, everything that follows could be repeated in situations where the Ricci bound only holds on a smaller set.
\end{remark}

We also recall the following lemma of Simon and Topping \cite[lemma 3.2]{SimonTopping22}.

\begin{lemma}\label{PROP shrinking balls}
    Let $(\mfd^n,g(t))_{t\in[0,T)}$ be a (possibly incomplete) Ricci flow. Suppose that for $t\in[0,T)$, we have
    \[
    \ric{g(t)} \leq \frac{3}{2}(n-1)f^2(t) \hspace{3ex}\text{in}\hspace{3ex} B_{g(0)}(p,r),
    \]
    or merely in $B_{g(0)}(p,r)\cap B_{g(t)}\left(p,r-4(n-1)F(t)\right)$, where $f:[0,T)\to[0,\infty)$ is a continuous integrable function, and $F(t):=\int_0^tf$. Then
    \[
    B_{g(t)}\left(p,r-4(n-1)F(t)\right) \setin B_{g(0)}(p,r)
    \hspace{8ex} \text{for all }t\in[0,T).
    \]
\end{lemma}


The two results above allow us to construct natural cutoff functions on $\mfd$.

\begin{lemma}\label{PROP cutoff func}
    Let $(\mfd^n,g(t))_{t\in[0,T)}$ be a (possibly incomplete) Ricci flow, with $B_{g(0)}(p,r)\Subset\mfd$ for some $p\in\mfd$, $r>0$. Let $f:[0,T)\to[0,\infty)$ be a continuous integrable function, and $F(t):=\int_0^tf$. Suppose that on $\mfd\times[0,T)$, we have
    \[
    \ric{g(t)} \leq \frac{3}{2}(n-1)f^2(t).
    \]
    Let $0<\rho<r$, and assume (by reducing $T$ if necessary) that $r-\rho-4(n-1)F(t)>1/f(t)$ for all $t\in[0,T)$.

    Given some constant $A>0$, and some weakly decreasing function $\sigma:[0,1]\longrightarrow[0,\infty)$, let $\twiddle{\psi}\in C^1(\reals)$ be such that the following holds.
    \begin{equation}\label{EQN conditions on psi twiddle}
    \begin{aligned}
        &\twiddle{\psi}\equiv0 &\text{in } (-\infty,0),\\
        &\twiddle{\psi}\equiv1 &\text{in } (1,\infty), \\
        &0\leq\twiddle{\psi}'\leq \twiddle{\psi}\,\sigma(\twiddle{\psi}) &\text{in } \reals, \\
        &\twiddle{\psi}''\geq-A\twiddle{\psi} &\text{in } \reals,
    \end{aligned}
    \end{equation}
    where the final statement holds in the barrier sense. That is, for each $s\in\reals$, there is a smooth function $\twiddle{\Psi}$ defined near $s$, such that $\twiddle{\Psi}\geq\twiddle{\psi}$ with equality at $s$, and $\twiddle{\Psi}''\geq-A\twiddle{\Psi}$ at $s$. Define $\psi_{\rho,r}:\mfd\times[0,T)\longrightarrow\reals$ by
    \[
    \psi_{\rho,r}(x,t) := e^{-A\rho^{-2}t} \ \twiddle{\psi}\left(\frac{r-d(x,t)-4(n-1)F(t)}{\rho}\right).
    \]

    Then $\psi_{\rho,r}$ enjoys the following properties.
    \begin{align*}
        \psi_{\rho,r}&\equiv e^{-A\rho^{-2}t} \hspace{8ex} &\text{in } B_{g(t)}(p,r-\rho-4(n-1)F(t)), \\
        \psi_{\rho,r}&\equiv 0   &\text{outside } B_{g(t)}(p,r-4(n-1)F(t)), \\
        |\nabla\psi_{\rho,r}| &\leq \rho^{-1} \psi_{\rho,r} \sigma\left(\psi_{\rho,r} \right) &\text{on }\mfd, \\
        \heatop\psi_{\rho,r}&\leq0 &\text{in the barrier sense},
    \end{align*}
    where by the penultimate inequality we mean that for any $(x,t)$ and any $C^2$ function $\Psi_{\rho,r}$ which satisfies $\Psi_{\rho,r}(x,t)=\psi_{\rho,r}(x,t)$ and $\Psi_{\rho,r}\leq\psi_{\rho,r}$ near $(x,t)$, we have $|\nabla\Psi_{\rho,r}| \leq \rho^{-1} \Psi_{\rho,r} \sigma\left(\Psi_{\rho,r} \right)$ at $(x,t)$.
\end{lemma}

\begin{proof}
    The first two properties follow by the construction of $\psi_{\rho,r}$.

    The gradient bound follows by considering the difference quotient of any suitable barrier function, and noting that $\sigma$ is weakly decreasing.

    For the final property, for fixed $(x,t)$ we take the barrier
    \[
    \Psi_{\rho,r}(x,t) := e^{-A\rho^{-2}t} \ \twiddle{\Psi}\left(\frac{r-D(x,t)-4(n-1)F(t)}{\rho}\right),
    \]
    where $\twiddle{\Psi}$ is the barrier for $\twiddle{\psi}$ given above, and $D$ is the barrier for $d$ given by Lemma \ref{PROP Perelman's lemma}. The result follows by a straightforward calculation, using Lemmas \ref{PROP Perelman's lemma} and \ref{PROP shrinking balls}.
\end{proof}


We now discuss some choices of the function $\twiddle{\psi}$ which we can use in Lemma \ref{PROP cutoff func}.

\begin{lemma}\label{PROP Specific cutoff functions}
    Let $m>0$. Then there exists $C_1>0$, and a function $\twiddle{\psi}^{*,m}\in C^1(\reals)$ such that \eqref{EQN conditions on psi twiddle} holds with
    \begin{equation}\label{EQN psi 1 properties}
    \begin{split}
        \sigma(s) &= mC_1s^{-1/m}\\
        A&=mC_1.
        \end{split}
    \end{equation}

    Further, given $m>0$, $\eps\in(0,1/2]$, there exists a function $\twiddle{\psi}^{\eps,m}\in C^1(\reals)$ such that \eqref{EQN conditions on psi twiddle} holds with
    \begin{equation}\label{EQN psi epsilon properties}
    \begin{split}
        \sigma(s) &= m\eps^{-1}\left(\max\left\{ |m^{-1}\log s|, 2^{2/\eps} \right\}\right)^{1+\eps},\\
        A &= 2^{-1+2/\eps}m.
    \end{split}
    \end{equation}
\end{lemma}

\begin{proof}
    First, suppose we already have some $\twiddle{\psi}$ which satisfies \eqref{EQN conditions on psi twiddle} for some $A$ and some function $\sigma$. Then where $\twiddle{\psi}>0$, we have
    \[
    \left|(\twiddle{\psi}^m)'\right|
    = m \twiddle{\psi}^{m-1}\left|\twiddle{\psi}'\right|
    \leq m \twiddle{\psi}^{m-1} \twiddle{\psi} \sigma(\twiddle{\psi})
    = \twiddle{\psi}^m \cdot m \sigma\left((\twiddle{\psi}^m)^{1/m}\right),
    \]
    and
    \[
    (\twiddle{\psi}^m)''
    = m \twiddle{\psi}^{m-1} \twiddle{\psi}'' + m(m-1) \twiddle{\psi}^{m-2} (\twiddle{\psi}')^2
    \geq -mA \twiddle{\psi}^m.
    \]
    So \eqref{EQN conditions on psi twiddle} holds for $\twiddle{\psi}^m$, with $\sigma(s)$ replaced with $m\sigma\left(s^{1/m}\right)$, and $A$ replaced with $mA$. Therefore, it is enough to settle the case $m=1$.

    The first part of the lemma now follows by finding a cutoff function with bounded derivative, and choosing $C_1$ appropriately.

    Now we turn our attention to the second part of the lemma, and build a set of more refined cutoff functions. Note that if we take $\twiddle{\psi}$ to satisfy $\twiddle{\psi}(s)=\exp\left(-s^{-1/\eps}\right)$ for small $s$, then (in this region at least), we have that $\sigma(s)=\frac{1}{\eps}|\log s|^{1+\eps}$.

    We now need to construct the rest of the cutoff function in such a way that $\sigma(s)$ is not too large away from the region where $\twiddle{\psi}$ is small, and that $A$ is well controlled.

    To this end, take constants $\delta_0,\delta_1,\delta_2\in(0,1)$, and $A>0$ with $A\delta_2<\pi/2$, all to be chosen later, and define $\twiddle{\psi}$ as follows (see Figure \ref{PIC cutoff function}).

\[
\twiddle{\psi}(s)=
\begin{cases}
    0 & -\infty<s\leq0 \\
    \exp\left(-s^{-1/\eps}\right) & 0<s\leq\delta_0 \\
    \exp\left( -\delta_0^{-1/\eps} + \eps^{-1}\delta_0^{-1-1/\eps}(s-\delta_0) \right) & \delta_0<s\leq\delta_0+\delta_1 \\
    \cos\left( A(\delta_0+\delta_1+\delta_2-s) \right) & \delta_0+\delta_1<s\leq\delta_0+\delta_1+\delta_2 \\
    1 & \delta_0+\delta_1+\delta_2<s<\infty
\end{cases}
\]

\begin{figure}[ht]
   \centering
   \def\svgscale{1}
   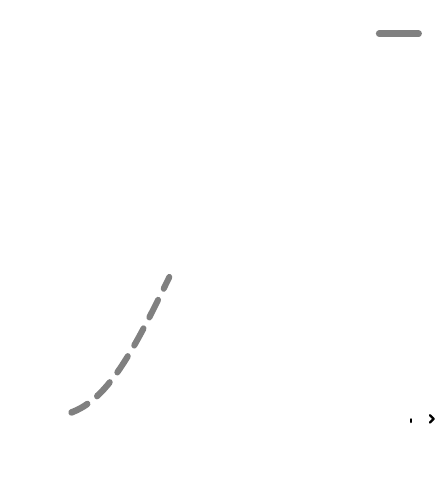
   \caption{A picture of $\twiddle{\psi}$}
    \label{PIC cutoff function}
\end{figure}

Then provided $\twiddle{\psi}$ is $C^1$, and $\delta_0+\delta_1+\delta_2\leq1$, $\twiddle{\psi}$ is a valid choice of function for Lemma \ref{PROP cutoff func}, and \eqref{EQN conditions on psi twiddle} holds with

\[
\sigma(s) = \eps^{-1}\big(\max\big\{ |\log s|, \delta_0^{-1/\eps} \big\}\big)^{1+\eps}
\]
and our choice of $A$ above.

We now use the requirement that $\twiddle{\psi}$ is $C^1$ to determine some of the constants in terms of $\delta_0$. First note that $\twiddle{\psi}$ is $C^1$ at $s=\delta_0$. Then for $\twiddle{\psi}$ to be $C^1$ at $s=\delta_0+\delta_1$, we must have
\begin{align*}
    \cos\left( A\delta_2 \right)
    &=\exp\left( \delta_0^{-1/\eps} \left( \frac{\delta_1}{\eps\delta_0}-1 \right) \right),\\
    A\sin\left( A\delta_2 \right)
    &=\eps^{-1}\delta_0^{-1-1/\eps} \exp\left( \delta_0^{-1/\eps} \left( \frac{\delta_1}{\eps\delta_0}-1 \right) \right).
\end{align*}
We choose $\delta_1$, $\delta_2$ in terms of $\delta_0$, $A$ such that these both hold. Since $\delta_0^{-1/\eps}>1$, this certainly means that $\delta_1<\eps\delta_0$. Recall also that $A\delta_2<\pi/2$. Therefore, to be sure of meeting the restriction $\delta_0+\delta_1+\delta_2\leq1$, we can simply choose $\delta_0$, $A$ such that
\[
(1+\eps)\delta_0 + \frac{\pi}{2A} = 1.
\]

We are left with some flexibility in the final choice. It will be advantageous for the proof of Theorem \ref{PROP main theorem 1/t^1+gamma and 1} to choose $A,\delta_0$ such that
\[
A=\half\delta_0^{-1/\eps},
\]
so we take this value for $A$, and we take $\delta_0$ such that
\[
(1+\eps)\delta_0 + \pi\delta_0^{1/\eps} = 1.
\]
Since $\eps\in(0,1/2]$, we have $\delta_0\geq2^{-2}$, and the result follows.
\end{proof}

\section{Existing results}\label{SEC existing results}
We now introduce $\ph$, some scalar function which is initially non-negative, and for which we want to find an upper bound. 

The following calculation is a generalisation of a calculation in \cite[Proof of Theorem 1.1]{LeeTam22}. We will need to repeatedly use it later on, so we isolate it as a lemma. Here $\eta$ is to be viewed as the candidate upper bound for $\ph$.

\begin{lemma} \label{PROP faffy calculation}
    Let $(\mfd^n,g(t))_{t\in[0,T)}$ be a (possibly incomplete) Ricci flow, with $B_{g(0)}(p,r)\Subset\mfd$ for some $p\in\mfd$, $r>0$. Let $f:[0,T)\to[0,\infty)$ be a continuous integrable function, and $F(t):=\int_0^tf$. Suppose that on $\mfd\times[0,T)$, we have
    \[
    \ric{g(t)} \leq \frac{3}{2}(n-1)f^2(t).
    \]
    Let $0<\rho<r$, and assume (by reducing $T$ if necessary) that $r-\rho-4(n-1)F(t)>1/f(t)$ for all $t\in[0,T)$.
    
    Suppose we also have continuous functions $\ph,L:\mfd\times[0,T)\to\reals$, such that on $B_{g(t)}(p,r-4(n-1)F(t))$ we have
    \[
    \heatop\ph\leq L\ph
    \]
    wherever $\ph>0$, in the barrier sense.
    
    Given some $C^1$ function $\eta:[0,T)\to(0,\infty)$, and a cutoff function $\psi_{\rho,r}$ as given by Lemma \ref{PROP cutoff func}, write $u:=-\psi_{\rho,r}\ph+\eta$. Assume that $u>0$ for $t\in(0,\eps)$, where the given $\eps>0$ is not required to be independent of $\eta$, or any other quantity.

    Then if $t_0\in(0,T]$ is the first time for which $u=0$ somewhere, and $x_0\in\mfd$ is a point where $u(x_0,t_0)=0$, then
    \begin{equation}\label{EQN faffy calc result}
        \left( \log \eta \right)' \leq L + \frac{2}{\rho^2} \sigma^2\left( \frac{\eta}{\ph} \right)
    \end{equation}
    at $(x_0,t_0)$, where $\sigma$ is as in Lemma \ref{PROP cutoff func}.
\end{lemma}

\begin{proof}
    At $(x_0,t_0)$, we can use the barriers for $\psi_{\rho,r}$, $\ph$ to construct a barrier $U$ for $u$. Then $\nabla U=0$, $\heatop U\leq0$ at $(x_0,t_0)$. The result follows from a brief calculation.
\end{proof}

The following theorem is the beginning of several steps of refinement, which will ultimately give our main result. In the case $\rho=\frac{r}{2}$, it is proven along the way to \cite[Theorem 1.1]{LeeTam22}, and the proof for general $\rho>0$ is very similar, so we omit it. Note that the cutoff function used in this case is as given by Lemma \ref{PROP cutoff func} using $\twiddle{\psi}^{*,1}$.

\begin{theorem}\label{PROP main theorem constant bound}
    Let $(\mfd^n,g(t))_{t\in[0,T)}$ be a (possibly incomplete) Ricci flow, with $B_{g(0)}(p,r)\Subset\mfd$ for some $p\in\mfd$, $r>0$. Let $f:[0,T)\to[0,\infty)$ be a continuous integrable function, and $F(t):=\int_0^tf$. Suppose that on $\mfd\times[0,T)$, we have
    \[
    \ric{g(t)} \leq \frac{3}{2}(n-1)f^2(t).
    \]

    Suppose we also have continuous functions $\ph,L:\mfd\times[0,T)\to\reals$, and $\alpha,\beta\geq0$ such that
    \[
    L\leq\alpha,
    \]
    \[
    \ph\leq\beta,
    \]
    \[
    \heatop\ph\leq L\ph \ \text{ wherever }\ph>0, \text{ in the barrier sense}
    \]
    on $B_{g(t)}(p,r-4(n-1)F(t))$ for all $t\in[0,T)$. Finally, suppose that $\ph\leq0$ on $B_{g(0)}(p,r)$ when $t=0$.

    Let $0<\rho<r$. Then for all $k>1$, there exists $\widehat{T}=\widehat{T}(f,r,\rho,\alpha,\beta,k)>0$ such that for $t\in[0,T)\cap[0,\widehat{T})$, we have
    \[
    \ph\leq t^k
    \hspace{8ex} \text{in } B_{g(t)}(p,r-\rho-4(n-1)F(t)).
    \]
\end{theorem}

The next result is an improvement over \cite[Theorem 1.1]{LeeTam22}, although a similar result 
can be obtained using the method in \cite{LeeTam19}.

\begin{theorem}\label{PROP main theorem 1/t bound}
    Let $(\mfd^n,g(t))_{t\in[0,T)}$ be a (possibly incomplete) Ricci flow, with $B_{g(0)}(p,r)\Subset\mfd$ for some $p\in\mfd$, $r>0$. Let $f:[0,T)\to[0,\infty)$ be a continuous integrable function, and $F(t):=\int_0^tf$. Suppose that on $\mfd\times[0,T)$, we have
    \[
    \ric{g(t)} \leq \frac{3}{2}(n-1)f^2(t).
    \]

    Suppose we also have continuous functions $\ph,L:\mfd\times[0,T)\to\reals$, and $\alpha,\beta\geq0$ such that
    \[
    L\leq\frac{\alpha}{t},
    \]
    \[
    \ph\leq\frac{\beta}{t},
    \]
    \[
    \heatop\ph\leq L\ph \ \text{ wherever }\ph>0, \text{ in the barrier sense}
    \]
    on $B_{g(t)}(p,r-4(n-1)F(t))$ for all $t\in[0,T)$. Finally, suppose that $\ph\leq0$ on $B_{g(0)}(p,r)$ when $t=0$.

    Let $0<\rho<r$. Then there is a universal constant $C$, and some time $\widehat{T}=\widehat{T}(f,r,\rho,\alpha,\beta)>0$ such that for $t\in[0,T)\cap[0,\widehat{T})$, we have
    \[
    \ph\leq\exp\left( -\frac{C\rho^2}{t\log(1/t)} \right)
    \hspace{8ex} \text{in } B_{g(t)}(p,r-\rho-4(n-1)F(t)).
    \]
\end{theorem}

\begin{proof}
    We begin by largely replicating the final part of the proof of \cite[Theorem 1.1]{LeeTam22}, so we sketch this part. We fix $k\geq1$ for now, and seek a bound of $t^k$. Since $L$, $\ph$ are continuous back to time $0$, we may apply Theorem \ref{PROP main theorem constant bound}  to see that $\ph\leq t^{k+1}$ in $B_{g(t)}(p,r-\rho/2-4(n-1)F(t))$, for some short time which we have no control over.

    Take the function $\twiddle{\psi}^{*,4km}$ from Lemma \ref{PROP Specific cutoff functions}, where $m\geq1$ is to be chosen later. Use Lemma \ref{PROP cutoff func} to produce a cutoff function $\psi:=\psi_{\rho/2,\,r-\rho/2}$, which is zero outside $B_{g(t)}(p,r-\rho/2-4(n-1)F(t))$ and equal to $\exp(-Ckm\rho^{-2}t)$ in $B_{g(t)}(p,r-\rho-4(n-1)F(t))$, for some constant $C>0$. Recall that $\sigma(s)=Ckms^{-1/(2km)}$ in this case.

    After possibly reducing $T$, we may apply Lemma \ref{PROP faffy calculation} with candidate bound $\eta=2^{-4km}t^k$ and the cutoff function used above. Using the bounds on $L,\ph$, we see that at the first time $t_0>0$ where $\psi\ph+\eta=0$, we have
    \[
    \frac{k}{t_0} \leq \frac{\alpha}{t_0} + \frac{Ck^2m^2}{\rho^2} \left(\frac{2^{-4km}t_0^{k+1}}{\beta}\right)^{-1/(2km)}
    \]
    for some universal constant $C$. Assuming that $k$ is sufficiently large (depending on $\alpha$, $\beta$), this implies
    \[
    \frac{C\rho^2}{m^2k} \leq t_0^{1-1/m}.
    \]
    
    Then reducing $T$ if needed so that $\psi\geq1/2$ in $B_{g(t)}(p,r-\rho-4(n-1)F(t))$, we have shown that
    \[
    \ph\leq t^k \hspace{5ex} \text{in } B_{g(t)}(p,r-\rho-4(n-1)F(t)), \hspace{5ex} \text{provided } t^{1-1/m}\leq\frac{C\rho^2}{m^2k}.
    \]
    This is true for any $k\geq K(\alpha,\beta), m\geq1$.

    Now, we optimise over $k,m$ to obtain the best possible bound for $\ph$. Fix a time $t\in(0,\min\{T,1\})$. For any fixed $m\geq1$, the best bound is obtained by maximising $k$, which requires us to take
    \[
    k=\frac{C\rho^2}{m^2t^{1-1/m}},
    \]
    provided $t$ is small enough that $k\geq K(\alpha,\beta)$. Then, it turns out that to maximise this choice of $k$ over $m$, we must take
    \[
    m=\half\log(1/t),
    \]
    provided $t$ is also sufficiently small that $m\geq1$. Then
    \[
    k=\frac{C\rho^2}{t\left(\log(1/t)\right)^2}.
    \]
    So we see that in $B_{g(t)}(p,r-\rho-4(n-1)F(t))$, we have
    \[
    \ph\leq t^k
    =t^{\,\wedge}\left(\frac{C\rho^2}{t\left(\log(1/t)\right)^2}\right)
    =\exp\left(-\frac{C\rho^2}{t\log(1/t)}\right),
    \]
    provided $t$ is sufficiently small, depending on $\alpha,\beta$.
\end{proof}

\section{Main result}\label{SEC main result}

We now make use of the more refined cutoff functions obtained in \S\ref{SEC cutoff functions} to prove a version of Theorem \ref{PROP main theorem 1/t bound} with the sharp bound. For now we assume a constant bound on $\ph$, which will be relaxed later.

\begin{theorem}\label{PROP main theorem 1/t^1+gamma and 1}
    Let $(\mfd^n,g(t))_{t\in[0,T)}$ be a (possibly incomplete) Ricci flow, with $B_{g(0)}(p,r)\Subset\mfd$ for some $p\in\mfd$, $r>0$. Let $f:[0,T)\to[0,\infty)$ be a continuous integrable function, and $F(t):=\int_0^tf$. Suppose that on $\mfd\times[0,T)$, we have
    \[
    \ric{g(t)} \leq \frac{3}{2}(n-1)f^2(t).
    \]

    Suppose we also have continuous functions $\ph,L:\mfd\times[0,T)\to\reals$, and constants $\alpha\geq0$, $\gamma\in\left[1/2,1\right)$ such that
    \[
    L\leq\alpha t^{-(1+\gamma)},
    \]
    \[
    \ph\leq1,
    \]
    \[
    \heatop\ph\leq L\ph \ \text{ wherever }\ph>0, \text{ in the barrier sense}
    \]
    on $B_{g(t)}(p,r-4(n-1)F(t))$ for all $t\in[0,T)$. Finally, suppose that $\ph\leq0$ on $B_{g(0)}(p,r)$ when $t=0$.

    Let $0<\rho<r$. Then there is a constant $C=C(\gamma)$, and some time $\widehat{T}=\widehat{T}(f,r,\rho)>0$ such that
    \[
    \ph \leq \exp\left(-C\rho^2t^{-1}\right) \hspace{8ex} \text{in } B_{g(t)}(p,r-\rho-4(n-1)F(t))
    \]
    for $t\in[0,T)\cap[0,\widehat{T}]\cap\left[0,(\half C\rho^2\alpha^{-1}\gamma)^{1/(1-\gamma)}\right]$.
\end{theorem}

\begin{remark}
    In particular, $\widehat{T}$ is chosen such that $r-\rho-4(n-1)F(t)>1/f(t)$ for all $t\in[0,\widehat{T})$.
\end{remark}

\begin{remark} \label{RMK Can do exp t^-gamma bound}
    The proof below also works if we assume a bound of $\exp(\beta t^{-\gamma})$ on $\ph$, but this makes it harder to estimate the constant $C(\gamma)$, and we only require the above result in order to continue.
\end{remark}


\begin{proof}
    The proof proceeds in two steps - first, we essentially repeat the method from \cite{LeeTam22} to obtain a bound for $\ph$ of the form $\exp\left(kt^{-\gamma}\right)$, where the time for which the bound is valid depends on $k$. Then, for fixed $t>0$ we optimise over $k$ and any other constants which appear, in order to obtain our final bound.

    \textbf{Step 1 - Method from \cite{LeeTam22}}\\ 
    First, since $L$, $\ph$ are continuous back to $t=0$, given any $\delta\in(0,1)$, we can apply Theorem \ref{PROP main theorem 1/t bound} to see that
    \[
    \ph\leq\exp\left(-\frac{C\delta^2\rho^2}{t\log(1/t)}\right) \hspace{8ex}\text{in }B_{g(t)}(p,r-\delta\rho-4(n-1)F(t))
    \]
    for some short time which we have no control over.
    
    Now take the function $\twiddle{\psi}^{\eps,m}$ from Lemma \ref{PROP Specific cutoff functions}, for $\eps\in(0,1/2]$, $m>0$ to be chosen later. Use Lemma \ref{PROP cutoff func} to produce a cutoff function $\psi:=\psi_{(1-\delta)\rho,\,r-\delta\rho}$, which is zero outside $B_{g(t)}(p,r-\delta\rho-4(n-1)F(t))$ and equal to $\exp(-2^{-1+2/\eps}m(1-\delta)^{-2}\rho^{-2}t)$ in $B_{g(t)}(p,r-\rho-4(n-1)F(t))$. Recall that $\sigma(s)=m\eps^{-1}\left(\max\left\{ |m^{-1}\log s|, 2^{2/\eps} \right\}\right)^{1+\eps}$ in this case.

    Choose $T_1$ such that $r-\rho-4(n-1)F(t)>\frac{1}{f(t)}$ for $t\in[0,T_1]$, and write $\widehat{T}_\delta:=\min\left\{T_1,(1-\delta)^2\rho^2\right\}$. Then we also have $\exp\left((1-\delta)^{-2}\rho^{-2}t\right)\geq e^{-1}$ for $t\in[0,\widehat{T}_\delta]$,
    so $\psi\geq\exp\left(-2^{-1+2/\eps}m\right)$ in $B_{g(t)}(p,r-\rho-4(n-1)F(t))$ for $t\in[0,\widehat{T}_\delta]$.

    Now we apply Lemma \ref{PROP faffy calculation} with candidate bound $\eta(t):=\exp\left(-2^{-1+2/\eps}m-kt^{-\gamma}\right)$, for $k>0$ to be chosen later, and the cutoff function above. Using the bounds on $L,\ph$, we see that at the first time $t_0>0$ where $\psi\ph+\eta=0$, we have
    \[
    (\gamma k-\alpha)t_0^{-(1+\gamma)} 
    \leq 
    \frac{2m^2\eps^{-2}}{(1-\delta)^2\rho^2} 
    \left(\max\left\{ m^{-1}kt_0^{-\gamma}+2^{-1+2/\eps},
    2^{2/\eps} \right\} \right)^{2+2\eps}.
    \]
    Now if $t_0$ is small enough that $m^{-1}kt_0^{-\gamma}>2^{-1+2/\eps}$, then $m^{-1}kt_0^{-\gamma}+2^{-1+2/\eps}>2^{2/\eps}$, so
    \begin{align*}
        (\gamma k-\alpha)t_0^{-(1+\gamma)} &< 2(1-\delta)^{-2}\rho^{-2}  m^2\eps^{-2}\left(2m^{-1}kt_0^{-\gamma}\right)^{2+2\eps} \\
        &= 2^{3+2\eps}(1-\delta)^{-2}\rho^{-2}  m^{-2\eps}\eps^{-2}\left(kt_0^{-\gamma}\right)^{2+2\eps}.
    \end{align*}
    Now further assume that $\gamma k\geq2\alpha$, so that
    \[
    \gamma t_0^{-1} < 2^{4+2\eps}(1-\delta)^{-2}\rho^{-2}  m^{-2\eps}\eps^{-2}\left(kt_0^{-\gamma}\right)^{1+2\eps}.
    \]
    Finally, choose $\eps\in(0,1/2]$ such that $1+2\eps=\gamma^{-1}$, but continue to also write $\eps$ for now. Then
    \[
    \gamma < 2^{4+2\eps}(1-\delta)^{-2}\rho^{-2}  m^{-2\eps}\eps^{-2}k^{1+2\eps}.
    \]
    if the first time $t_0$ when $\psi\ph=\eta$ satisfies $kt_0^{-\gamma}>2^{-1+2/\eps}m$, and $\gamma k\geq2\alpha$.
    
    For fixed $k\geq2\alpha\gamma^{-1}$, $m>0$, take $\tau_1>0$ such that $k\tau_1^{-\gamma}= 2^{-1+2/\eps}m$, and suppose $m$, $k$ were chosen such that
    \[
    \gamma \geq 2^{4+2\eps}(1-\delta)^{-2}\rho^{-2}  m^{-2\eps}\eps^{-2}k^{1+2\eps}.
    \]
    If $t\in(0,\tau_1)$, then $kt^{-\gamma}>2^{-1+2/\eps}m$, so $t$ cannot be the first time $t_0$ when $\psi\ph=\eta$.
    Recall that $\psi\geq\exp\left(-2^{-1+2/\eps}m\right)$ in $B_{g(t)}(p,r-\rho-4(n-1)F(t))$ for $t\in[0,\widehat{T}_\delta]$. 
    We have therefore shown that $\ph\leq\exp\left(-kt^{-\gamma}\right)$ in $B_{g(t)}(p,r-\rho-4(n-1)F(t))$ for $t\in[0,T)\cap[0,\widehat{T}_\delta]\cap[0,\tau_1]$. 
    Sending $\delta\to0$ gives the following.
    
    \begin{claim}\label{PROP Claim of t^-gamma bound}
        Write $\eps=(\gamma^{-1}-1)/2$. Then for any $m>0$, $k\geq2\alpha\gamma^{-1}$, such that
        \begin{equation}\label{EQN restriction on k, m}
            \gamma \geq 2^{4+2\eps}\rho^{-2}  m^{-2\eps}\eps^{-2}k^{1+2\eps},
        \end{equation}
        the following holds.
        
        Choose $T_1$ such that $r-\rho-4(n-1)F(t)>\frac{1}{f(t)}$ for $t\in[0,T_1)$, and write $\widehat{T}:=\min\left\{T_1,\rho^2\right\}$. Take $\tau_1>0$ such that $k\tau_1^{-\gamma}=2^{-1+2/\eps}m$. Then
        \[
        \ph\leq\exp\left(-kt^{-\gamma}\right) \hspace{8ex} \text{in } B_{g(t)}(p,r-\rho-4(n-1)F(t))
        \]
        for $t\in[0,T)\cap[0,\widehat{T}]\cap[0,\tau_1]$.
    \end{claim}

    \textbf{Step 2 - Optimisation}\\ 
    Now we use the result of Claim \ref{PROP Claim of t^-gamma bound}, and optimise over $m$, $k$ to obtain an $\exp(-Ct^{-1})$ bound.
    
    Fix some $t\in[0,T)\cap[0,\widehat{T}]$. For fixed $k\geq2\alpha\gamma^{-1}$, choose $m>0$ such that $kt^{-\gamma}=2^{-1+2/\eps}m$. Then provided \eqref{EQN restriction on k, m} holds, we have
    \[
    \ph(\cdot,t)\leq\exp\left(-kt^{-\gamma}\right) \hspace{8ex} \text{in } B_{g(t)}(p,r-\rho-4(n-1)F(t)).
    \]
    
    Recalling that $2\eps=\gamma^{-1}-1$, we have $k^{2\eps}\left(t^{-\gamma}\right)^{-1+1/\gamma}=\left(2^{-1+2/\eps}m\right)^{2\eps}$, so $2^{2\eps}m^{-2\eps}k^{2\eps}=2^4t^{1-\gamma}$, and we can rewrite \eqref{EQN restriction on k, m} as
    \begin{equation}\label{EQN restriction on k}
        2^{-8}\rho^2\eps^2\gamma \geq kt^{1-\gamma}.
    \end{equation}
    We may take $k$ to give equality in \eqref{EQN restriction on k} provided $t$ is small enough that
    \[
    2\alpha\gamma^{-1}\leq k= 2^{-8}\rho^2\eps^2\gamma t^{\gamma-1}.
    \]

    Therefore, taking
    \[
    C:=C(\gamma)=2^{-8}\eps^2\gamma=2^{-10}(\gamma^{-1}-1)^2\gamma,
    \]
    we have
    \[
    \ph(\cdot,t)\leq\exp\left(-kt^{-\gamma}\right)=\exp\left(-C\rho^2t^{-1}\right)  \hspace{8ex} \text{in } B_{g(t)}(p,r-\rho-4(n-1)F(t)),
    \]
    provided
    \[
    2\alpha\gamma^{-1}\leq C\rho^2 t^{\gamma-1}.
    \]
    The result follows.
\end{proof}

Theorem \ref{PROP main theorem 1/t^1+gamma and 1} gives an $\exp(-Ct^{-1})$ bound, in the case where we have the bounds $L\leq\alpha t^{-(1+\gamma)}$, $\ph\leq\exp(\beta t^{-\gamma})$, for $\gamma$ strictly less than $1$ (in principle - see Remark \ref{RMK Can do exp t^-gamma bound}). Our next task is to improve to the case $\gamma=1$. It turns out that this can be done, but only if $\alpha$, $\beta$ are small enough.

\begin{theorem}\label{PROP main theorem 1/t^2 and exp -1/t}
    There is a universal constant $\widehat{\beta}$ such that the following holds.

    Let $(\mfd^n,g(t))_{t\in[0,T)}$ be a (possibly incomplete) Ricci flow, with $B_{g(0)}(p,r)\Subset\mfd$ for some $p\in\mfd$, $r>0$. Let $f:[0,T)\to[0,\infty)$ be a continuous integrable function, and $F(t):=\int_0^tf$. Suppose that on $\mfd\times[0,T)$, we have
    \[
    \ric{g(t)} \leq \frac{3}{2}(n-1)f^2(t).
    \]

    Suppose we also have continuous functions $\ph,L:\mfd\times[0,T)\to\reals$, and a constant $\beta\leq\widehat{\beta}$, such that
    \[
    L\leq\big(\widehat{\beta}-\beta\big)\rho^2 t^{-2},
    \]
    \[
    \ph\leq\exp\left(\beta\rho^2 t^{-1}\right),
    \]
    \[
    \heatop\ph\leq L\ph \ \text{ wherever }\ph>0, \text{ in the barrier sense}
    \]
    on $B_{g(t)}(p,r-4(n-1)F(t))$ for all $t\in[0,T)$. Finally, suppose that $\ph\leq0$ on $B_{g(0)}(p,r)$ when $t=0$.

    Let $0<\rho<r$. Then there is some time $\widehat{T}=\widehat{T}(f,r,\rho)>0$ such that
    \[
    \ph \leq \exp\big(-\big(\widehat{\beta}-\beta\big)\rho^2t^{-1}\big) \hspace{8ex} \text{in } B_{g(t)}(p,r-\rho-4(n-1)F(t))
    \]
    for $t\in[0,T)\cap[0,\widehat{T}]$.
\end{theorem}

\begin{remark}
    In particular, $\widehat{T}$ is as given by Theorem \ref{PROP main theorem 1/t^1+gamma and 1}, that is it is chosen such that $r-\rho-4(n-1)F(t)>1/f(t)$ for all $t\in[0,\widehat{T})$.
\end{remark}

\begin{proof}
    Before we begin, we define a useful function. Take $\zeta:[0,\infty)\to[0,\infty)$ to be the $C^1$ function
    \[
    \zeta_1(t):=
    \begin{cases}
        \frac{1}{4}(x+1)^2&0\leq x<1,\\x&1\leq x<\infty
    \end{cases}
    \]
    (see Figure \ref{PIC zeta1}), and for any $t_0>0$, take $\zeta_{t_0}(t):=t_0\zeta\left(\frac{t}{t_0}\right)$, so that
    \[
    t\leq \zeta_{t_0}(t)\leq\max\left\{t_0,t\right\},
    \]
    \[
    0\leq \zeta_{t_0}'(t)\leq1,
    \]
    \[
    \min\left\{\frac{1}{t_0^2},\frac{1}{t^2}\right\} \leq \frac{\zeta_{t_0}'(t)}{\zeta_{t_0}^2(t)} \leq \frac{1}{t^2}
    \]
    for all $t\geq0$.

    \begin{figure}[ht]
        \centering
        \def\svgscale{1}
        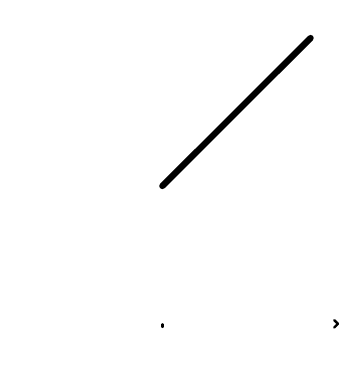
        \caption{A graph of $\zeta_1$}
        \label{PIC zeta1}
    \end{figure}
    
    We begin by establishing the following.
    \begin{claim}\label{PROP Claim of transfer to exp t^-gamma}
        Suppose the assumptions of Theorem \ref{PROP main theorem 1/t^2 and exp -1/t} hold, but we instead have $\ph\leq\exp(\beta\rho^2 t^{-\gamma})$ for some $\beta>0,\gamma\in[1/2,1)$, and $L=0$. Then we have that
        \[
        \ph \leq \exp\left(\beta\rho^2 t^{-\gamma}-C\rho^2t^{-1}\right) \hspace{8ex} \text{in } B_{g(t)}(p,r-\rho-4(n-1)F(t))
        \]
        for $t\in[0,T)\cap[0,\widehat{T}]\cap\left[0,(\half C\beta^{-1})^{1/(1-\gamma)}\right]$, where $C=C(\gamma)$ is as in Theorem \ref{PROP main theorem 1/t^1+gamma and 1}.
    \end{claim}
    \begin{proof}[Proof of Claim \ref{PROP Claim of transfer to exp t^-gamma}.]
        Since $\ph$ is continuous back to $t=0$, there must be some $t_0>0$ such that $\ph\leq\exp(\beta\rho^2 t_0^{-\gamma})$ for all $(x,t)\in S$. Then consider $\twiddle{\ph}:=\exp\left(-\beta\rho^2 \zeta_{t_0}^{-\gamma}\right)\ph$. We have
        \[
        \twiddle{\ph} \leq 
        \exp\left(-\beta\rho^2 (\max\{t_0,t\})^{-\gamma}\right) \cdot
        \min\left\{\exp(\beta\rho^2 t_0^{-\gamma}),\exp(\beta\rho^2 t^{-\gamma})\right\}
        =1
        \]
        in $S$, and
        \[
        \heatop\twiddle{\ph} \leq \gamma\beta\rho^2 \zeta_{t_0}^{-(1+\gamma)}\zeta_{t_0}'\twiddle{\ph}.
        \]
        Note that $\gamma\beta\rho^2 \zeta_{t_0}^{-(1+\gamma)}\zeta_{t_0}'$ is continuous on $[0,T)$, with
        \[
        \gamma\beta\rho^2\zeta_{t_0}^{-(1+\gamma)}\zeta_{t_0}' \leq \gamma\beta\rho^2 t^{-(1+\gamma)}.
        \]
        
        Hence we can apply Theorem \ref{PROP main theorem 1/t^1+gamma and 1} with $\alpha=\gamma\beta\rho^2$ to see that
        \[
        \twiddle{\ph} \leq \exp\left(-C\rho^2t^{-1}\right) \hspace{8ex} \text{in } B_{g(t)}(p,r-\rho-4(n-1)F(t))
        \]
        for $t\in[0,T)\cap[0,\widehat{T}]\cap\left[0,(\half C\beta^{-1})^{1/(1-\gamma)}\right]$.

        Therefore $\ph\leq\exp\left(\beta\rho^2 \zeta_{t_0}^{-\gamma}-C\rho^2t^{-1}\right)\leq\exp\left(\beta\rho^2 t^{-\gamma}-C\rho^2t^{-1}\right)$ for $(x,t)$ in the same set.
    \end{proof}

    Now we choose $\widehat{\beta}$, and use Claim \ref{PROP Claim of transfer to exp t^-gamma} to establish the result in the case $\alpha=0$, $\beta=\widehat{\beta}$.

    \begin{claim}\label{PROP Claim of transfer to exp t^-1}
        There is a universal constant $\widehat{\beta}$ such that the following holds.

        Suppose the assumptions of Theorem \ref{PROP main theorem 1/t^2 and exp -1/t} hold with $L=0$. Then
        \[
        \ph \leq 1 \hspace{8ex} \text{in } B_{g(t)}(p,r-\rho-4(n-1)F(t))
        \]
        for $t\in[0,T)\cap[0,\widehat{T}]$.
    \end{claim}

    \begin{proof}[Proof of Claim \ref{PROP Claim of transfer to exp t^-1}.]
        Fix the choice of $\gamma\in[1/2,1)$ which gives the largest possible value of $C(\gamma)$ in Theorem \ref{PROP main theorem 1/t^1+gamma and 1}, and take $\widehat{\beta}:=\half C(\gamma)$. Write 
        \[
        S:=\{(x,t):t\in[0,T)\cap[0,\widehat{T}],x\in B_{g(t)}(p,r-4(n-1)F(t))\},
        \]
        where $\widehat{T}$ is from Theorem \ref{PROP main theorem 1/t^1+gamma and 1}.
        
        By continuity, $\ph$ must remain bounded as $t\searrow0$, so either $\ph$ attains its maximum over $(x,t)\in S$ at some time $t_0$, or $\ph\leq\exp\big(\widehat{\beta}\rho^2T^{-1}\big)$ for all $(x,t)\in S$, in which case we take $t_0:=T$.
    
        Then for $(x,t)\in S$, we have
        \[
        \ph \leq \min \big\{ \exp\big(\widehat{\beta}\rho^2t_0^{-1}\big), \exp\big(\widehat{\beta}\rho^2t^{-1}\big) \big\} \leq \exp\big(\widehat{\beta}\rho^2t_0^{-(1-\gamma)}t^{-\gamma}\big).
        \]
        So applying Claim \ref{PROP Claim of transfer to exp t^-gamma} with $\beta=\widehat{\beta}t_0^{-(1-\gamma)}$ gives
        \[
        \ph \leq \exp\big(\widehat{\beta}\rho^2t_0^{-(1-\gamma)} t^{-\gamma}-C\rho^2t^{-1}\big) \hspace{8ex} \text{in } B_{g(t)}(p,r-\rho-4(n-1)F(t))
        \]
        for $t\in[0,T)\cap[0,\widehat{T}]\cap\left[0,\big(\half C\widehat{\beta}^{-1}\big)^{1/(1-\gamma)}t_0\right]$.
        
        Recall that $C=C(\gamma)=2\widehat{\beta}$. If $t_0=T$, then $\ph\leq1$ for all $(x,t)\in S$ as required. Otherwise, we have $\ph\leq1$ when $t=t_0$. Since by assumption $\ph$ attains its maximum at this time, the result follows.
    \end{proof}

    Finally, we use Claim \ref{PROP Claim of transfer to exp t^-1} to establish Theorem \ref{PROP main theorem 1/t^2 and exp -1/t}. Take $L$, $\ph$ satisfying the conditions of Theorem \ref{PROP main theorem 1/t^2 and exp -1/t}. Then by continuity, there is some time $t_0\in[0,T)\cap[0,\widehat{T}]$ such that $\ph\leq\exp\left(\beta\rho^2 t_0^{-1}\right)$, $L\leq\big(\widehat{\beta}-\beta\big)\rho^2 t_0^{-2}$. This means that if $\tau\in(0,t_0]$ then
    \[
    \ph \leq \exp\left(\beta\rho^2 \zeta_\tau^{-1}\right), \hspace{8ex} L \leq \big(\widehat{\beta}-\beta\big)\rho^2 \zeta_\tau'\zeta_\tau^{-2},
    \]
    for all $(x,t)\in S$.

    Let $\twiddle{\ph}:=\exp\big(\big(\widehat{\beta}-\beta\big)\rho^2\zeta_\tau^{-1} \big)\ph$. Then we still have $\twiddle{\ph}=0$ initially, and $\ph\leq\exp\big(\widehat{\beta}\rho^2\zeta_\tau^{-1} \big)\leq\exp\big(\widehat{\beta}\rho^2t^{-1} \big)$ for all $(x,t)\in S$. Finally,
    \[
    \heatop \twiddle{\ph} \leq -(\widehat{\beta}-\beta)\rho^2\zeta_\tau'\zeta_\tau^{-2} \twiddle{\ph} + \big(\widehat{\beta}-\beta\big)\rho^2\zeta_\tau'\zeta_\tau^{-2}\twiddle{\ph} \leq 0
    \]
    in $B_{g(t)}(p,r-4(n-1)F(t))$ for all $t\in [0,T)$. Applying Claim \ref{PROP Claim of transfer to exp t^-1} to $\twiddle{\ph}\leq1$ gives
    \[
    \ph\leq\exp(-\big(\widehat{\beta}-\beta\big)\rho^2\zeta_\tau^{-1}) \hspace{8ex} \text{in } B_{g(t)}(p,r-\rho-4(n-1)F(t))
    \]
    for $t\in[0,T)\cap[0,\widehat{T}]$.
    Finally, sending $\tau\to0$ gives the result.
\end{proof}

Theorem \ref{PROP main thm for powers of t} now follows quickly.

\begin{proof}[Proof of Theorem \ref{PROP main thm for powers of t}.]
    Apply Theorem \ref{PROP main theorem 1/t^2 and exp -1/t} with $\frac{3}{2}(n-1)f^2(t)=\alpha t^{-\gamma}$ and $\rho=r/2$, allowing the value of $\widehat{\beta}$ to decrease as needed.
\end{proof}


\printbibliography

\end{document}

%% file: CutoffFunction.pdf_tex
\begingroup%
  \makeatletter%
  \providecommand\color[2][]{%
    \errmessage{(Inkscape) Color is used for the text in Inkscape, but the package 'color.sty' is not loaded}%
    \renewcommand\color[2][]{}%
  }%
  \providecommand\transparent[1]{%
    \errmessage{(Inkscape) Transparency is used (non-zero) for the text in Inkscape, but the package 'transparent.sty' is not loaded}%
    \renewcommand\transparent[1]{}%
  }%
  \providecommand\rotatebox[2]{#2}%
  \newcommand*\fsize{\dimexpr\f@size pt\relax}%
  \newcommand*\lineheight[1]{\fontsize{\fsize}{#1\fsize}\selectfont}%
  \ifx\svgwidth\undefined%
    \setlength{\unitlength}{213.63178884bp}%
    \ifx\svgscale\undefined%
      \relax%
    \else%
      \setlength{\unitlength}{\unitlength * \real{\svgscale}}%
    \fi%
  \else%
    \setlength{\unitlength}{\svgwidth}%
  \fi%
  \global\let\svgwidth\undefined%
  \global\let\svgscale\undefined%
  \makeatother%
  \begin{picture}(1,1.08377554)%
    \lineheight{1}%
    \setlength\tabcolsep{0pt}%
    \put(0,0){\includegraphics[width=\unitlength,page=1]{CutoffFunction.pdf}}%
    \put(0.92313492,0.0842484){\makebox(0,0)[t]{\lineheight{1.25}\smash{\begin{tabular}[t]{c}1\end{tabular}}}}%
    \put(0,0){\includegraphics[width=\unitlength,page=2]{CutoffFunction.pdf}}%
    \put(0.01882902,1.00374781){\makebox(0,0)[t]{\lineheight{1.25}\smash{\begin{tabular}[t]{c}1\end{tabular}}}}%
    \put(0,0){\includegraphics[width=\unitlength,page=3]{CutoffFunction.pdf}}%
    \put(0.61597688,0.02276938){\makebox(0,0)[t]{\lineheight{1.25}\smash{\begin{tabular}[t]{c}$\delta_2$\end{tabular}}}}%
    \put(0,0){\includegraphics[width=\unitlength,page=4]{CutoffFunction.pdf}}%
    \put(0.27029571,0.02276938){\makebox(0,0)[t]{\lineheight{1.25}\smash{\begin{tabular}[t]{c}$\delta_1$\end{tabular}}}}%
    \put(0,0){\includegraphics[width=\unitlength,page=5]{CutoffFunction.pdf}}%
    \put(0.10970047,0.02276938){\makebox(0,0)[t]{\lineheight{1.25}\smash{\begin{tabular}[t]{c}$\delta_0$\end{tabular}}}}%
    \put(0.02132854,0.0842726){\makebox(0,0)[t]{\lineheight{1.25}\smash{\begin{tabular}[t]{c}0\end{tabular}}}}%
    \put(0,0){\includegraphics[width=\unitlength,page=6]{CutoffFunction.pdf}}%
  \end{picture}%
\endgroup%

%% file: zeta.pdf_tex
\begingroup%
  \makeatletter%
  \providecommand\color[2][]{%
    \errmessage{(Inkscape) Color is used for the text in Inkscape, but the package 'color.sty' is not loaded}%
    \renewcommand\color[2][]{}%
  }%
  \providecommand\transparent[1]{%
    \errmessage{(Inkscape) Transparency is used (non-zero) for the text in Inkscape, but the package 'transparent.sty' is not loaded}%
    \renewcommand\transparent[1]{}%
  }%
  \providecommand\rotatebox[2]{#2}%
  \newcommand*\fsize{\dimexpr\f@size pt\relax}%
  \newcommand*\lineheight[1]{\fontsize{\fsize}{#1\fsize}\selectfont}%
  \ifx\svgwidth\undefined%
    \setlength{\unitlength}{167.04052302bp}%
    \ifx\svgscale\undefined%
      \relax%
    \else%
      \setlength{\unitlength}{\unitlength * \real{\svgscale}}%
    \fi%
  \else%
    \setlength{\unitlength}{\svgwidth}%
  \fi%
  \global\let\svgwidth\undefined%
  \global\let\svgscale\undefined%
  \makeatother%
  \begin{picture}(1,1.05365275)%
    \lineheight{1}%
    \setlength\tabcolsep{0pt}%
    \put(0,0){\includegraphics[width=\unitlength,page=1]{zeta.pdf}}%
    \put(0.46604216,0.04847673){\makebox(0,0)[t]{\lineheight{1.25}\smash{\begin{tabular}[t]{c}1\end{tabular}}}}%
    \put(0,0){\includegraphics[width=\unitlength,page=2]{zeta.pdf}}%
    \put(0.02408084,0.51128076){\makebox(0,0)[t]{\lineheight{1.25}\smash{\begin{tabular}[t]{c}1\end{tabular}}}}%
    \put(0.02727754,0.04850767){\makebox(0,0)[t]{\lineheight{1.25}\smash{\begin{tabular}[t]{c}0\end{tabular}}}}%
    \put(0,0){\includegraphics[width=\unitlength,page=3]{zeta.pdf}}%
  \end{picture}%
\endgroup%